\documentclass[11pt, a4paper]{amsart}
\usepackage{amsmath,amssymb,amscd,amsfonts}

\newtheorem{thm}{Theorem}[section]

\newtheorem{rem}[thm]{Remarks}
\newtheorem{prop}[thm]{Proposition}

\newcommand{\Pic}{{\text{\rm Pic}}}

\newcommand{\OO}{{\mathcal {O}}}

\newcommand{\inv}{^{-1}}

\newcommand{\Ga}{\Gamma}

\newcommand{\epsi}{\epsilon}
\newcommand{\fie}{\varphi}
\newcommand{\si}{\sigma}
\newcommand{\Si}{\Sigma}

\newcommand{\Q}{{\mathbb Q}}
\newcommand{\Z}{{\mathbb Z}}
\newcommand{\F}{{\mathbb F}}
\newcommand{\pp}{{\mathbb{P}}}

\newcommand{\tpi}{\tilde{\pi}}
\newcommand{\tsi}{\tilde{\si}}
\numberwithin{equation}{section}

\title{The bicanonical map of surfaces with  $p_g=0$ and $K^2\ge 7$, II}
\author{Margarida Mendes Lopes \and Rita Pardini}
\date{}
\begin{document}

\begin{abstract} We study the  minimal complex  surfaces of general type with $p_g=0$ and
$K^2=7$ or $8$ whose bicanonical map is not birational. We show that if  $S$ is such a surface,
then the bicanonical map has degree 2 (see \cite{mp}) and there is a fibration
$f\colon S\to \pp^1$ such that: i) the general fibre $F$ of $f$ is a genus $3$ hyperelliptic curve;
ii) the involution induced by the bicanonical map of $S\/$ restricts to the
hyperelliptic involution of $F$.\par
Furthermore, if $K^2_S=8$, then  $f$ is an isotrivial fibration with   6 double fibres,
and  if $K^2_S=7$, then $f$ has 5 double fibres  and it has precisely  one fibre
with reducible support, consisting of two components.

\noindent 2000 Mathematics Classification: 14J29.
\end{abstract}
\maketitle
\section{Introduction}

 A minimal surface $S$ of general type with
$p_g(S)=0$
 satisfies the inequalities $1\le K^2_S\le 9$.  It is known that for $K_S^2\geq
2$ the image of the bicanonical map $\fie$ of $S$ is a surface and that for
$K^2_S\geq 5$  the bicanonical map is always a morphism. In \cite{mp} it is shown  that
 $\fie$ is birational if  $K_S^2=9$ and that
 the degree of $\fie$ is at most $2$  if $K_S^2=7$ or
$K_S^2=8$.\par In this note, which is a complement to \cite{mp}, we describe the
surfaces  with $K^2_S=7$ or $8$ for which the degree of $\fie$ is  $2$. We prove the
following:

\begin{thm}\label{main} Let $S$ be a minimal smooth complex surface of general
type with
$p_g(S)=0$ and
$K_S^2=7$ or $8$ for which the bicanonical map is not birational. Then:
\begin{itemize}
\item[i)]  $K_S$ is ample;
\item[ii)]  there is a fibration $f\colon S\to\pp^1$ such that
 the general fibre $F$ of $f$ is hyperelliptic of genus 3;
\item[iii)] the bicanonical involution of $S$  induces the hyperelliptic involution on
$F$.
\end{itemize}
\smallskip

\noindent Furthermore
\begin{itemize}

\item[iv)] if $K^2_S=8$, then  $f$ is an isotrivial fibration with   6 double fibres;
\item[v)] if $K^2_S=7$, then $f$ has 5 double fibres  and it has precisely  one fibre
with reducible support, consisting of two components.
\end{itemize}
\end{thm}

We remark that surfaces satisfying the assumptions of Theorem \ref{main} do exist. In
\cite{mp} we have given examples of such surfaces  with $K^2_S=7$ and with $K^2_S=8$, and
the case $K^2_S=8$ is completely classified in \cite{rita}.
\medskip

\noindent{\bf Acknowledgements.}  The present collaboration takes place in the framework
of the european contract EAGER, no. HPRN-CT-2000-00099. The first author is a member of
CMAF and of the Departamento de Matem\'atica da Faculdade de Ci\^encias da Universidade de
Lisboa and the second author is a member of GNSAGA of CNR.

\medskip
\noindent{\bf Notation and conventions.}  We work over the complex numbers;
all
varieties are assumed to be  compact and algebraic.  We do not distinguish
between line bundles and divisors on a smooth variety, using the additive
and the multiplicative notation interchangeably. Linear equivalence is denoted by
$\equiv$. A {\em nodal curve} or {\em $-2-$curve} $C$  on a surface is a curve
isomorphic to
$\pp^1$ and such that $C^2=-2$.  The rest of  the notation is standard in algebraic
geometry.

\section{Involutions and the bicanonical map}

In this section we introduce some notation   and we give  a
 criterion  to decide whether    an involution $\si$ of a surface $S$ of
general type with $p_g=0$ acts  trivially on   the bicanonical system  of $S$ (cf.
Proposition \ref{fixnumber}).

Let
$S$ be a smooth surface. An {\em involution} of
$S$ is an automorphism
$\si$ of
$S$ of order 2. The fixed locus of $\si$ is the union  of a smooth curve $R$ and of $k$
isolated points $P_1\dots P_k$. We denote by $\pi\colon S\to  \Si:=S/\si$ the projection
onto the quotient,  by
$B$ the image of $R$ and by $Q_i$ the image of $P_i$, $i=1\dots k$. The surface $\Si$ is
normal and $Q_1\dots Q_k$   are ordinary double points, which are the only singularities
of $\Si$. In particular, the singularities of $\Si$ are canonical and the adjunction
formula gives
$K_S=\pi^*K_{\Si}+R$. Let $\epsi\colon X\to S$ be the blow-up of $S$ at $P_1\dots P_k$
and let
$E_i$ be the exceptional curve over $P_i$, $i=1\dots k$. It is easy to check that $\si$
induces an involution
$\tsi$ of
$X$ whose fixed locus is the union of $R_0:=\epsi\inv R$ and of $E_1\dots E_k$. Denote by
$\tpi\colon X\to Y:=X/\tsi$ the projection onto the quotient and set $B_0:=\tpi(R_0)$,
$C_i:=\tpi(E_i)$, $i=1\dots k$. The surface $Y$ is smooth and the $C_i$ are
disjoint $-2-$curves. Denote  by $\eta\colon Y\to
\Si$ the morphism induced by
$\epsi$. The map $\eta$ is the minimal resolution of the singularities of $\Si$ and
there is  a commutative diagram:
\begin{equation}\label{diagram}
\begin{CD}\ X@>\epsi>>S\\ @V\tpi VV  @VV\pi V\\ Y@>\eta >> \Si
\end{CD}
\end{equation}
\bigskip

 We say that a map $\psi\colon S\to W$ is composed with $\si$ if $\psi\circ \si=\psi$.
\begin{prop}\label{fixnumber} Let $S$ be a minimal surface of general type with
$p_g(S)=0$ and let
$\si$ be an involution of $S$. Assume that $|2K_S|$ has no fixed
component   and let
$\fie\colon S\to \pp^{K^2_S}$ be the bicanonical map.

Then $\fie$  is composed with  $\si$ if and only if
   the number of isolated fixed points of $\si$ is equal to  $K^2_S+4$.
\end{prop}
\begin{proof}

We use the notation  introduced above for involutions.

\noindent We have
$\chi(\OO_S)>0$, since $S$ is of general type, and
$p_g(S)=0$ by assumption, hence
$q(S)=0$ and, as a consequence,
$p_g(\Si)=q(\Si)=0$. We also have $P_2(S)=\chi(\OO_S)+K^2_S=1+K^2_S$.

If   $B=0$,   then $K_S=\pi^*K_{\Si}$, hence   $K_{\Si}$ is nef and big and $\Si$ is
a surface of general type. The standard formulas for double covers give:
$$1=\chi(\OO_S) =2\chi(\OO_\Si)-k/4=2-k/4,$$ i.e. $k=4<K^2_S+4$. In addition we have
$$P_2(S)=K^2_S+1=2K^2_{\Si}+1>K^2_{\Si}+1=P_2(\Si),$$ hence $\fie$ is not composed with
$\si$.

Now   assume that $B\ne 0$ and consider diagram (\ref{diagram}).  By the theory of double
covers,  there exists
$L\in
\Pic(Y)$ such that
$2L\equiv B_0+\sum C_i$ and
$\tilde{\pi}_*\OO_X=\OO_Y\oplus L\inv$. The adjunction formula gives
$$2K_X=\tilde{\pi}^*(2K_Y+B_0+\sum C_i)$$ and the projection formulas for double covers
give
$$H^0(X,2K_X)= H^0(Y,2K_Y+L) \oplus H^0(Y,2K_Y +B_0+\sum C_i).$$ The bicanonical map
$\fie$  is composed with $\si$ iff either $H^0(Y,2K_Y+L)=0$ or $H^0(Y,2K_Y +B_0+\sum
C_i)=0$. The elements of $H^0(Y,2K_Y+L)$ pull-back on $S$ to the ``odd'' sections of
$2K_X$, which vanish on $R_0$. On the other hand,
$$|2K_X|=\epsi^*|2K_S|+2\sum E_i,$$  hence  the fixed part of
$|2K_X|$ is supported on
$\sum E_i$, since   by assumption $|2K_S|$   has no base
component. It follows that
$H^0(Y,2K_Y+B_0+\sum C_i)\ne 0$, and $\fie$ is composed with $\si$ iff
$H^0(Y,2K_Y+L)=0$.

We remark that $2K_{\Si}+B$ is nef and big, since $2K_S=\pi^*(2K_{\Si}+B)$ is nef and big
by assumption.  We have  the equality of
$\Q-$divisors:
$$K_Y+L=\frac{1}{2}(2K_Y+B)+\frac{1}{2}\sum C_i.$$ The divisor
$\frac{1}{2}(2K_Y+B_0)=\frac{1}{2}\eta^*(2K_{\Si}+B)$ is nef and big and the divisor
$\frac{1}{2}\sum C_i$ is effective with normal crossings support and zero integral part.
Thus  $h^i(Y,2K_Y+L)=0$ for $i>0$ by  Kawamata--Viehweg vanishing and so:

\begin{equation}\label{1}  h^0(Y,2K_Y+L)=\chi(2K_Y+L)
=1+K^2_Y+\frac{3}{2}K_YL+\frac{1}{2}L^2.
\end{equation} Using again the projection formulas for double covers, we get
$$1=\chi(\OO_X)=\chi(\OO_Y)+\chi(K_Y+L)=1+\chi(K_Y+L),$$ which by Riemann--Roch is
equivalent to
\begin{equation}\label{2} L^2+LK_Y=-2.
\end{equation} Finally, we have: $$k=K^2_S-K^2_X=K^2_S-2(K_Y+L)^2.$$ Using equalities
(\ref{1}) and (\ref{2}),  this equation can be rewritten: $$k=K^2_S+4- 2h^0(2K_Y+L).$$
Summing up, $\fie$ is composed with $\si$ iff $h^0(Y,2K_Y+L)=0$ iff $k=K^2_S+4$.
\end{proof}
\section{The main result} Throughout all this section we assume that $S$ is a minimal
surface of general type with
$p_g(S)=0$ and $K^2_S=7$ or $8$ and that the bicanonical map $\fie\colon
S\to\pp^{K^2_S}$ is not birational.

By Theorem 1 of \cite{mp}, the degree of $\fie$  is equal to 2.  We keep the notation of
the previous section. In particular, we denote by
$\si$ the involution  associated with $\fie$ and  by $\pi\colon S\to \Si:=S/\si$ the  quotient
map.
\begin{prop}\label{fibref}
 $\Si$ is a rational surface with $K^2_{\Si}=-4$ whose singularities  are $K_S^{2}+4$
ordinary double points.  Furthermore:
\begin{itemize}
    \item[i)] there is a fibration $h\colon \Si\to\pp^1$  whose general fibre  is smooth
rational;
\item[ii)] $h$ has $6$ double fibres for $K^2_S=8$ and it has $5$ double  fibres for
$K^2_S=7$;
\item [iii)] every double fibre of $h$ contains at least 2 singular points;
\item[iv)] if $K^{2}_S=8$, then all the fibres of $h$ have irreducible  support. If
$K^2_S=7$, then  $h$ has exactly one fibre with  reducible support, consisting of two
components.
\end{itemize}
\begin{proof} By Theorem 3 of \cite{xiaodeg}, $\Si$ is a rational surface. By Reider's
theorem, the bicanonical system of $S$ has no base points, hence, by Proposition
\ref{fixnumber},
  the singularities of   $\Si$ are $K_{S}^{2}+4$  ordinary double
    points. The rank
$\nu$ of
$NS(\Si)$ is
$\le
    h^{2}(S)=10-K^{2}_S$, hence we have $\nu\le 2$ if $K^2_S=8$ and $\nu\le 3$ if
$K^2_S=7$ and  we can apply the results of
\cite{nodes}. By  Lemma 3.1  and Theorem 3.2  of \cite{nodes}, there
    exists a fibration $h\colon \Si\to \pp^{1}$ such that the general fibre of $h$ is a
smooth rational curve and such that
     $h$ has
$6$ double fibres for $K^{2}_S=8$ and   it has $5$ double fibres for $K^{2}_S=7$.
Furthermore Theorem 3.3 of
    \cite{nodes} implies that $\nu\ge 2$ if $K^2_S=8$ and that $\nu\ge 3$ if
    $K^2_S=7$. Hence we have $\nu=2$ if $K^2_S=8$,  $\nu=3$ if $K^2_S=7$
    and
    $K^{2}_{\Si}=-4$ in both cases.  Let $\eta \colon Y\to\Si$ be the minimal
    resolution and consider the fibration  $\tilde{h}:=h\circ\eta\colon
    Y\to\pp^{1}$. It
    is possible to factor $\tilde{h}$
    as $Y\to Y^{\prime}\to\pp^{1}$, where $Y\to Y^{\prime}$ is a birational morphism of
smooth surfaces and
$Y^{\prime}\to
\pp^{1}$ is a relatively
    minimal fibration. It follows that $Y^{\prime}$ is a ruled surface $\F_{e}$ and
    $Y\to Y^{\prime}$ is a sequence of 12 blow-ups.  On the other hand, $\eta$ contracts
12 vertical
    curves for $K^{2}_S=8$ and 11 vertical curves for $K^{2}_S=7$. Therefore
     all the fibres of $h$ have irreducible support in the former
    case and   there is precisely one  fibre with reducible support, consisting of two
    components, in the latter case.
    \end{proof}

\end{prop}
\begin{thm}\label{fibreg} There is a fibration $f\colon S\to\pp^1$ such that:
\begin{itemize}
\item[i)] the general fibre $F$ of $f$ is hyperelliptic of genus 3;
\item[ii)] $\si$ induces the hyperelliptic involution on $F$;
\item[iii)] if $K^2_S=8$, then  $f$ is an isotrivial fibration  whose singular
fibres are   6 double fibres;
\item[iv)] if $K^2_S=7$, then the multiple fibres of $f$ are  5 double fibres  and there
is   precisely  one fibre of $f$ with reducible support, consisting of two components.
\end{itemize}
\end{thm}
\begin{proof}  We set $f:=h\circ\pi$, where $h$ is the fibration of Proposition
\ref{fibref}.  Hence the general fibre $F$   of $f$ is hyperelliptic and $\si$ induces
the hyperelliptic involution on it. By Proposition \ref{fibref}, $f$ has 6 double fibres
if $K^2_S=8$ and it has $5$ double
 fibres if $K^2_S=7$.

Assume $K^2_S=8$ and let $\psi\colon C\to\pp^1$ be the double cover branched over the  6
image points of the double fibres of $f$. Taking fibre product and  normalization, one
gets a commutative diagram:
\begin{equation}\label{diagram2}
\begin{CD}\ X@>\tilde{\psi}>>S\\ @V\tilde{f} VV  @VV f V\\ C@>\psi>> \pp^1
\end{CD}
\end{equation} The map $\tilde{\psi}$ is \'etale, hence $X$ is smooth and we have:
$$\chi(\OO_X)=2\chi(\OO_S)=2,\quad
K^2_X=2K^2_S=16.$$  The fibrations $f$ and $\tilde{f}$ have the same general fibre, which
we
 still denote by $F$. Notice   that the genus $g(F)$  of $F$ is odd, since $f$ has double
fibres. By
\cite{appendice},
 we have:
\begin{equation}\label{3} 2=\chi(\OO_X)\ge (g(C)-1)(g(F)-1)=(g(F)-1),
\end{equation}
hence  $g(F)=3$ and (\ref{3}) is an equality. So, again by \cite{appendice},
  $\tilde{f}$ is smooth and isotrivial, hence $f$ is also
 isotrivial. Since $f$ is smooth, the singular fibres of $f$ occur only at the branch
points of
$\psi$, hence they are precisely the 6 double fibres.

Assume now that $K^2_S=7$. Then there exists a $\Z_2^2-$cover $\psi\colon  C\to \pp^1$
branched on the 5 image points of the double fibres of $f$ (cf. \cite{ritaabel},
Proposition 2.1). Taking base change and normalization, we obtain a commutative diagram
similar to (\ref{diagram2}). In this case $\tilde{\psi}$ is an \'etale
$\Z_2^2-$cover, hence $K^2_X=28$ and $\chi(\OO_X)=4$. Again \cite{appendice} gives
$$28=K^2_X\ge 8(g(C)-1)(g(F)-1)=8(g(F)-1),$$ namely
$g(F)=3$, since $g(F)$ is odd and $S$ is of general type. Assume now that the fibration
$f$ has another  multiple fibre $F_0$,  different from the 5 double fibres that we have
considered. Then $F_0$ is also a double fibre, since $K_SF_0=4$, and we can consider a
diagram similar to (\ref{diagram2}) where $\psi\colon C\to\pp^1$ is the double cover
branched on the points corresponding to these 6 double fibres. Again, $X$ is smooth and
\cite{appendice} gives
$$14=K^2_X\ge 8(g(C)-1)(g(F)-1)=16,$$
a contradiction. Hence for $K^2_S=7$ the multiple fibres of $f$ are precisely the 5
double fibres.

In addition,  we have   $h^2(S,\Z)=10-K^2_S=3$ and $f$ has at least a
reducible fibre $F_0$ by Proposition \ref{fibref}, iv).  Since $K_S$ and the components of
$F_0$ are independent in the N\'eron--Severi group, it follows that $F_0$ is the only
reducible fibre and that it has precisely two components.
\end{proof}

To complete the proof of Theorem \ref{main}, we need to show:
\begin{thm}\label{ample} The canonical class $K_S$ is ample.
\end{thm}
\begin{proof} If $K^2_S=8$, then the statement is an immediate consequence  of Miyaoka's
inequality (\cite{miyaoka}, \S 2).

 We assume,  from now on,  that $K^2_S=7$. We prove the theorem for this case, by a step-by-step
analysis which will show that $S$ does not contain $-2-$curves.
\medskip

\noindent Step 1: {\it
\smallskip If $C\subset S$ is a $-2-$curve, then it is contained in a fibre of $f$ (cf.
Theorem \ref{fibreg}).}

 Assume by contradiction $FC>0$, where $F$ denotes the fibre of $f$, as usual. By
Proposition
\ref{fibreg},
$f$ has  double fibres, hence
 $FC$ is even, say  $FC=2d$.  Then the Hurwitz formula  applied to the  map $f|_C$ gives
$-2\ge -2(2d)+5d=d$, a contradiction.
\medskip

Write $F_i=2D_i$, $i=1\dots 5$, for the double  fibres of $f$ and set
$H:=D_1+\dots +D_5+F$.  The divisor $K_S+H$ restricts to $K_F$ on  every fibre of
$f$ different from one of the   $F_i$ and it restricts to $K_{D_i}$ on
$D_i$ for $i=1\dots 5$.
\smallskip

\noindent Step 2: {\it The restriction maps $H^0(K_S+H)\to  H^0(K_F)$ and
$H^0(K_S+H)\to H^0(K_{D_i})$, $i=1\dots 5$, are surjective.}

 Let $F$ be a  fibre of
$f$ different from the $F_i$ and let
$\Delta\in |H|$ be the   divisor $F+D_1+\dots+D_5$. Consider the restriction
sequence:

$$0\to\OO_S(-H)\to \OO_S\to \OO_{\Delta}\to 0.$$
Since $S$ is regular, the corresponding long exact sequence in cohomology gives:

$$h^1(-H)=h^0(\OO_{\Delta})-1=h^0(\OO_F)+h^0(\OO_{D_1})+\dots +h^0(\OO_{D_5})-1=5,$$
where the last equality follows from the fact that the divisors $F$ and $D_1\dots D_5$
are
$1-$connected  by Theorem
\ref{fibreg} and by Zariski's Lemma (cf. \cite{bpv}, Ch. 3, \S 8). The same
argument shows:

$$h^1(-H+F)=h^1(-H+D_i)=4, \quad  i=1\dots 5.$$ The cokernel of the
restriction $H^0(K_S+H)\to  H^0(K_F)$  is the kernel of the map $H^1(K_S+H-F)\to
H^1(K_S+H)$. In turn, the cokernel of the latter map is the $1-$dimensional space
$H^1(K_F)$, since by Serre duality $h^2(K_S+H-F)=h^0(-(D_1+\dots +D_5))=0$. It follows
that
$H^1(K_S+H-F)\to H^1(K_S+H)$ is injective, since $h^1(K_S+H)=h^1(-H)=5$ and
$h^1(K_S+H-F)=h^1(-H+F)=4$. This shows that $H^0(K_S+H)\to  H^0(K_F)$  is onto. The
surjectivity of $H^0(K_S+H)\to H^0(K_{D_i})$, $i=1\dots 5$, can be proven in the same
way.
\medskip

\noindent Step 3: {\it The reducible fibre of $f$ (cf. Theorem \ref{fibreg}) contains
an effective  divisor
 $E$ with  $E^2=-1$, $K_SE=1$.}

\smallskip  A fibre of $f$ containing a divisor with  negative self-intersection is
necessarily reducible, hence  a divisor
$E$  as in the statement is contained in the only reducible fibre of $f$.

 Denote by
$\Phi$ the map given by $|K_S+H|$. Clearly $\Phi$ is composed with
$\si$, hence it has even degree. On the other hand, $(K_S+H)^2=35$ is odd and  thus
$|K_S+H|$ has base points.
 It follows that either  $|K_{D_i}|$ has base points for  some $i$ or there is a fibre
$F$, different from the $F_i$, such that
$|K_F|$ has base points.

Assume that, say, $|K_{D_1}|$ has base points. Then $D_1$ is not
$2-$connected (see, e.g., Proposition (A.7) of \cite{cfm}). Hence there is a
decomposition $D_1=A+B$, where $A$ and $B$  are effective and $AB\le 1$. On the other
hand  $D_1$  is $1-$connected,  hence we
have $AB=1$,
$A^2=B^2=-1$ and $K_SA=K_SB=1$, since $K_SD_1=2$.

Assume now that $|K_F|$ has base points, where $F$ is a fibre of $f$  different from
$F_1\dots F_5$. The same argument as above shows that $F$  decomposes as a sum   $A+B$,
where $A^2=B^2=-1$, $AB=1$ and $K_SA=1$,
$K_SB=3$.
\medskip

\noindent Step 4: {\it If $C\subset S$ is a $-2-$curve, then it is not contained in a
fibre of $f$.}

\smallskip Assume by contradiction that $C$ is contained in a fibre of $f$. Then $C$  is
a component of the only reducible fibre $F_0$ of $f$. On the other hand,  by Step 3,
$F_0$ contains an irreducible curve $\Ga$ with
$K_S\Ga=1$. By the adjunction formula,  one has either  $\Ga^2=-1$ or $\Ga^2=-3$.  By
Theorem \ref{fibreg}, iv), $F_0$ can be written as $4\Ga+mC$. An easy  analysis using
again Zariski's Lemma  shows that this is not  possible.
\medskip

The statement follows now from Steps 1 and 4.
\end{proof}
\begin{rem} i) In \cite{mp} we have presented examples of surfaces satisfying the
assumptions of  Theorem
\ref{fibreg}. In the example  with
$K_S^2=7$ the unique  fibre of the genus $3$ fibration having reducible support
is a double fibre. It would be interesting to know whether the reducible  fibre is
always one of the  double fibres.

 ii) The case $K^2_S=8$, where the fibration $f$ is
isotrivial, is completely classified in
\cite{rita}.

iii) All the known examples of surfaces satisfying the assumptions of Theorem \ref{fibreg}
are not isolated in the moduli space (see \cite{mp}, \cite{rita}), although the expected
dimension of the moduli space at the corresponding point is zero. The same phenomenon
occurs for Burniat surfaces, which have been characterized in \cite{mp2} as the only
minimal surfaces of general type with $K^2=6$ and $p_g=0$ such that the degree of the
bicanonical is 4, namely the largest possible degree for $3\le K^2\le 6$.\par

iv) It is interesting to remark that for $K^2\geq 6$, the bicanonical map having the
largest possible degree implies that the canonical class is ample (cf.
Theorem \ref{ample}, and
\cite{mp2}).

\end{rem}

\bigskip

\begin{tabbing} 1749-016 Lisboa, PORTUGALxxxxxxxxx\= 56127 Pisa, ITALY \kill Margarida
Mendes Lopes \> Rita Pardini\\ CMAF \> Dipartimento di Matematica\\
 Universidade de Lisboa \> Universit\a`a di Pisa \\ Av. Prof. Gama Pinto, 2 \> Via
Buonarroti 2\\ 1649-003 Lisboa, PORTUGAL \> 56127 Pisa, ITALY\\ mmlopes@lmc.fc.ul.pt \>
pardini@dm.unipi.it
\end{tabbing}
\end{document}